\documentclass[a4paper,11pt]{article}
\usepackage{amsmath}

\input{tcilatex}
\begin{document}

\title{Determining the Rolle function in Hermite interpolatory approximation
by solving an appropriate differential equation}
\author{J. S. C. Prentice \\
%EndAName
Senior Research Officer\\
Mathsophical Ltd.\\
Johannesburg, South Africa}
\maketitle

\begin{abstract}
We determine the pointwise error in Hermite interpolation by numerically
solving an appropriate differential equation, derived from the error term
itself. We use this knowledge to approximate the error term by means of a
polynomial, which is then added to the original Hermite polynomial to form a
more accurate approximation. An example demonstrates that improvements in
accuracy are significant.
\end{abstract}

\section{Introduction}

Recently, we reported on a technique for determining the Rolle function in
Lagrange interpolation, and how this could lead to an improvement in the
accuracy of the approximation \cite{Prentice 1}. In this short paper, we
extend that investigation to include Hermite interpolation. We consider the
same example as used in \cite{Prentice 1}, and show how significant
improvements in approximation accuracy can be achieved once the Rolle
function is known.

\section{Relevant Concepts}

Let $f\left( x\right) $ be a real-valued function. The \textit{Hermite
interpolating polynomial} $H_{2n+1}\left( x\right) $ of degree $2n+1,$ at
most, that interpolates the data $\{f\left( x_{0}\right) ,$ $f\left(
x_{1}\right) ,\ldots ,f\left( x_{n}\right) \}$ and $\{f^{\prime }\left(
x_{0}\right) ,$ $f^{\prime }\left( \ x_{1}\right) ,\ldots ,f^{\prime }\left(
x_{n}\right) \}$ at the nodes $\left\{ x_{0},x_{1},\ldots ,x_{n}\right\} ,$
where $x_{0}<x_{1}<\cdots <x_{n},$ has the properties%
\begin{align}
H_{2n+1}\left( x_{k}\right) & =f\left( x_{k}\right)  \label{interp prop f} \\
H_{2n+1}^{\prime }\left( x_{k}\right) & =f^{\prime }\left( x_{k}\right)
\label{interp prop f'}
\end{align}%
for $k=0,1,\ldots ,n.$ We have used the usual prime notation for
differentiation with respect to $x.$ We regard $H_{2n+1}\left( x\right) $ as
an approximation to $f\left( x\right) .$ The pointwise error in Hermite
interpolation, on $\left[ x_{0},x_{n}\right] ,$ is 
\begin{equation}
\Delta \left( x\left\vert H_{2n+1}\right. \right) \equiv f\left( x\right)
-H_{2n+1}\left( x\right) =\frac{f^{\left( 2n+2\right) }\left( \xi \left(
x\right) \right) }{\left( 2n+2\right) !}\dprod\limits_{k=0}^{n}\left(
x-x_{k}\right) ^{2},  \label{delta}
\end{equation}%
where $x_{0}<\xi \left( x\right) <x_{n},$ and may be derived by invoking
Rolle's Theorem \cite{I and K}\cite{K and C}. We necessarily assume here
that $f\left( x\right) $ is $\left( 2n+2\right) $-times differentiable. As
will be seen later, we must actually assume that $f\left( x\right) $ is $%
\left( 2n+3\right) $-times differentiable. We refer to $\xi \left( x\right) $
as the \textit{Rolle function}.

\section{The Rolle Function}

We employ the notation $Q_{n}\left( x\right) \equiv
\dprod\nolimits_{k=0}^{n}\left( x-x_{k}\right) $ and find, by
differentiating with respect to $x,$%
\begin{align*}
\left( 2n+2\right) !\left( f\left( x\right) -H_{2n+1}\left( x\right) \right)
& =f^{\left( 2n+2\right) }\left( \xi \left( x\right) \right) Q_{n}^{2}\left(
x\right) \\
\Rightarrow \left( 2n+2\right) !\left( f^{\prime }\left( x\right)
-H_{2n+1}^{\prime }\left( x\right) \right) & =2f^{\left( 2n+2\right) }\left(
\xi \right) Q_{n}Q_{n}^{\prime }\left( x\right) +Q_{n}^{2}\left( x\right) 
\frac{df^{\left( 2n+2\right) }\left( \xi \right) }{d\xi }\frac{d\xi }{dx} \\
& =2f^{\left( 2n+2\right) }\left( \xi \right) Q_{n}Q_{n}^{\prime }\left(
x\right) +Q_{n}^{2}\left( x\right) f^{\left( 2n+3\right) }\left( \xi \right) 
\frac{d\xi }{dx}.
\end{align*}%
In this expression, $f^{\left( 2n+2\right) }\left( \xi \right) $ denotes the 
$\left( 2n+2\right) $th derivative of $f\left( \xi \right) $ with respect to 
$\xi ,$ and similarly for $f^{\left( 2n+3\right) }\left( \xi \right) .$ We
now find%
\begin{equation*}
\frac{d\xi }{dx}=\frac{\left( 2n+2\right) !\left( f^{\prime }\left( x\right)
-H_{2n+1}^{\prime }\left( x\right) \right) -2f^{\left( 2n+2\right) }\left(
\xi \right) Q_{n}Q_{n}^{\prime }\left( x\right) }{Q_{n}^{2}\left( x\right)
f^{\left( 2n+3\right) }\left( \xi \right) }.
\end{equation*}%
If we have a particular value $\xi _{z}=\xi \left( x_{z}\right) $ available,
we have an initial-value problem that can be solved to yield the Rolle
function $\xi \left( x\right) .$ Note that the denominator in the above
expression requires the assumption that $f\left( x\right) $ is $\left(
2n+3\right) $-times differentiable.

\section{Numerical Example}

Consider the Hermite interpolation of%
\begin{eqnarray*}
f\left( x\right) &=&e^{x}\sin x \\
f^{\prime }\left( x\right) &=&e^{x}\sin x+e^{x}\cos x
\end{eqnarray*}%
over the nodes $\left\{ 0,\frac{3\pi }{2}\right\} .$ This is the same
example as used in \cite{Prentice 1}. Since $n=1$ we have 
\begin{equation*}
H_{3}\left( x\right) =ax^{3}+bx^{2}+c^{x}+d
\end{equation*}%
where the coefficients $a,b,c$ and $d$ are determined from the system%
\begin{equation*}
\left[ 
\begin{array}{cccc}
x_{0}^{3} & x_{0}^{2} & x_{0} & 1 \\ 
x_{1}^{3} & x_{1}^{2} & x_{1} & 1 \\ 
3x_{0}^{2} & 2x_{0} & 1 & 0 \\ 
3x_{1}^{2} & 2x_{1} & 1 & 0%
\end{array}%
\right] \left[ 
\begin{array}{c}
a \\ 
b \\ 
c \\ 
d%
\end{array}%
\right] =\left[ 
\begin{array}{c}
f\left( x_{0}\right) \\ 
f\left( x_{1}\right) \\ 
f^{\prime }\left( x_{0}\right) \\ 
f^{\prime }\left( x_{1}\right)%
\end{array}%
\right]
\end{equation*}%
with $x_{0}=0$ and $x_{1}=\frac{3\pi }{2}.$ We find $a=-2.8403,b=8.1595,c=1$
and $d=0$ (for ease of presentation, we quote numerical values to no more
than four decimal places, but all calculations were performed in double
precision).

Additionally,%
\begin{align*}
\Delta \left( x\left\vert H_{3}\right. \right) & =e^{x}\sin x-\left(
ax^{3}+bx^{2}+cx+d\right) \\
& =\frac{f^{\left( 4\right) }\left( \xi \left( x\right) \right) }{4!}\left(
x-x_{0}\right) ^{2}\left( x-x_{1}\right) ^{2} \\
& =-\frac{e^{\xi \left( x\right) }\sin \xi \left( x\right) }{6}\left(
x^{4}-3\pi x^{3}+\frac{9\pi ^{2}}{4}x^{2}\right)
\end{align*}%
so that%
\begin{equation}
\frac{d\xi }{dx}=\frac{18ax^{2}+12bx+6c-6e^{x}\left( \sin x+\cos x\right)
-A\left( x\right) e^{\xi }\sin \xi }{B\left( x\right) e^{\xi }\left( \sin
\xi +\cos \xi \right) }  \label{DE example 1}
\end{equation}%
where $A\left( x\right) \equiv 4x^{3}-9\pi x^{2}+\frac{9\pi ^{2}}{2}x$ and $%
B\left( x\right) \equiv x^{4}-3\pi x^{3}+\frac{9\pi ^{2}}{4}x^{2},$ and we
have used%
\begin{align*}
f^{\left( 4\right) }\left( \xi \right) & =-4e^{\xi }\sin \xi \\
f^{\left( 5\right) }\left( \xi \right) & =-4e^{\xi }\left( \sin \xi +\cos
\xi \right) .
\end{align*}

\medskip We solve this differential equation in a manner similar to that
used in \cite{Prentice 1}: we find an initial value at a point close to the
node $x_{0}=0$ (we cannot find $\xi _{z}$ at any interpolation node, because
the factor $\dprod\nolimits_{k=0}^{n}\left( x-x_{k}\right) ^{2}$ in (\ref%
{delta}) ensures that $\Delta \left( x_{z}\left\vert H_{2n+1}\right. \right)
=0$ at every interpolation node, \textit{regardless of the value of} $\xi ).$
Call this point $x_{z}$ and choose $x_{z}=10^{-5}.$ Since we know $f\left(
x\right) $ and $H_{3}\left( x\right) ,$ we can compute $\Delta \left(
x_{z}\left\vert H_{3}\right. \right) .$ Of course, this must be equal to 
\begin{equation*}
-\frac{e^{\xi _{z}}\sin \xi _{z}}{6}\left( x_{z}^{4}-3\pi x_{z}^{3}+\frac{%
9\pi ^{2}}{4}x_{z}^{2}\right)
\end{equation*}

\medskip where $\xi _{z}\equiv \xi \left( x_{z}\right) .$ We can easily
solve 
\begin{equation*}
\Delta \left( x_{z}\left\vert H_{3}\right. \right) =-\frac{e^{\xi _{z}}\sin
\xi _{z}}{6}\left( x_{z}^{4}-3\pi x_{z}^{3}+\frac{9\pi ^{2}}{4}%
x_{z}^{2}\right)
\end{equation*}%
numerically to find $\xi _{z}.$ In fact, we find two solutions $\xi
_{z}=0.9022$ and $\xi _{z}=3.0498.$ When we solve (\ref{DE example 1})
numerically, the first of these yields a Rolle function $\xi \left( x\right) 
$ that has negative values. This contradicts the constraint $x_{0}<\xi
\left( x\right) <x_{1},$ and so $\xi _{z}=0.9022$ is rejected as an initial
value. The second solution, on the other hand, gives an acceptable Rolle
function (see Figure 1). The numerical solution was obtained using a
seventh-order Runge-Kutta (RK) method \cite{Butcher} with a stepsize of $%
\sim 5\times 10^{-5},$ the same stepsize used in \cite{Prentice 1}.

Figure 2 shows the error curves - the LHS and RHS of (\ref{delta}) - for the
example. The curves are essentially indistinguishable. Figure 3 shows the
pointwise difference between these error curves. The difference is extremely
small, indicating the quality of our numerical solution of (\ref{DE example
1}), and the success of our algorithm for determining the Rolle function.

\section{Possible Applications}

Knowing the Rolle function $\xi \left( x\right) $ means we know $f^{\left(
2n+2\right) }\left( \xi \left( x\right) \right) .$ Hence, if we approximate $%
f^{\left( 2n+2\right) }\left( \xi \left( x\right) \right) $ by means of a
polynomial - perhaps a least-squares fit or a cubic spline - then, using (%
\ref{delta}), we have%
\begin{equation*}
f\left( x\right) \approx H_{2n+1}\left( x\right) +\frac{H_{\xi }\left(
x\right) }{\left( 2n+2\right) !}\dprod\limits_{k=0}^{n}\left( x-x_{k}\right)
^{2}\equiv H_{2n+1}\left( x\right) +E\left( x\right)
\end{equation*}%
where $H_{\xi }\left( x\right) $ denotes the polynomial that approximates $%
f^{\left( 2n+2\right) }\left( \xi \left( x\right) \right) ,$ and we have
implicitly defined the error polynomial $E\left( x\right) .$ The RHS of this
expression is simply a polynomial, and so constitutes a polynomial
approximation to $f\left( x\right) .$ Thus, our knowledge of $\xi \left(
x\right) $ allows us to improve the approximation $H_{2n+1}\left( x\right) $
by adding a polynomial term that approximates the pointwise error in $%
H_{2n+1}\left( x\right) .$

\subsection{The error polynomial}

For our earlier example, we have%
\begin{equation*}
E\left( x\right) =\frac{H_{\xi }\left( x\right) }{24}\left( x^{4}-3\pi x^{3}+%
\frac{9\pi ^{2}}{4}x^{2}\right) .
\end{equation*}%
We use the values of $\xi \left( x\right) $ from the RK process ($100000$
values over the interval $\left[ 0,\frac{3\pi }{2}\right] $) to generate $%
H_{\xi }\left( x\right) $ by fitting polynomials in a least-squares sense,
of varying degree. In Table 1, we show relevant results. The symbol $x_{i}$
denotes the RK nodes. The column "Max. error" shows 
\begin{equation*}
\max\nolimits_{i}\left\vert f\left( x_{i}\right) -\left( H_{3}\left(
x_{i}\right) +E\left( x_{i}\right) \right) \right\vert ,
\end{equation*}%
and $V$ is the variance of the fitted polynomial, 
\begin{equation*}
V\equiv \frac{\sqrt{\sum\nolimits_{i}\left( f^{\left( 4\right) }\left( \xi
\left( x_{i}\right) \right) -H_{\xi }\left( x_{i}\right) \right) ^{2}}}{%
100000}
\end{equation*}%
taken as a measure of goodness-of-fit.

\begin{center}
\begin{tabular}{|c|c|c|}
\hline
Degree of $H_{\xi }$ & Max. error & $V$ \\ \hline
$5$ & $9.6\times 10^{-3}$ & $2.1\times 10^{-5}$ \\ \hline
$7$ & $1.1\times 10^{-4}$ & $2.4\times 10^{-7}$ \\ \hline
$9$ & $3.0\times 10^{-6}$ & $9.6\times 10^{-9}$ \\ \hline
$11$ & $7.3\times 10^{-8}$ & $6.9\times 10^{-9}$ \\ \hline
\end{tabular}

Table 1: Relevant values pertaining to fitted polynomials.
\end{center}

\medskip

Clearly, the maximum approximation error decreases considerably as the
degree of $H_{\xi }$ increases. For reference, the maximum approximation
error for the original Hermite polynomial $H_{3}\left( x\right) $ is $7.04.$
We see that the use of $H_{\xi }$ improves the approximation by many orders
of magnitude. This effect was also observed in \cite{Prentice 1}. Note that
the degree of the error polynomial $E\left( x\right) $ is four plus the
degree of $H_{\xi }.$

We also consider the use of a cubic spline to generate $H_{\xi }.$ There are
several good reasons for this: we can use the RK values; the degree of $%
E\left( x\right) $ will be seven, at most; and, if we use a clamped spline,
we know the error bound in such an approximation \cite{B and F}\cite{Sch} is
given by 
\begin{equation*}
\frac{5\max_{i}\left\vert f^{\left( 8\right) }\left( x_{i}\right)
\right\vert }{384}h^{4}=1.14\times 10^{-16}
\end{equation*}%
where $h$ is the RK stepsize. In fact, we find 
\begin{equation*}
\max\nolimits_{i}\left\vert f\left( x_{i}\right) -\left( H_{3}\left(
x_{i}\right) +E\left( x_{i}\right) \right) \right\vert \sim 10^{-12}
\end{equation*}%
when using the cubic spline. We believe the discrepancy between this value
and the predicted bound is simply due to the less accurate values of $\xi
\left( x_{i}\right) $ generated by the RK method. This, of course, suggests
that the RK method could be a limiting factor in the overall accuracy of the
algorithm, and it would be appropriate to study how error control in said RK
method affects this accuracy. Not doing this here does not detract from our
demonstration, and so we will defer such a study to a future paper.

There is an important point to be made: 
\begin{align*}
H_{3}\left( x\right) +E\left( x\right) & =H_{3}\left( x\right) +\frac{H_{\xi
}\left( x\right) }{24}B\left( x\right) \\
H_{3}^{\prime }\left( x\right) +E^{\prime }\left( x\right) & =H_{3}^{\prime
}\left( x\right) +\frac{H_{\xi }^{\prime }\left( x\right) }{24}B\left(
x\right) +\frac{H_{\xi }\left( x\right) }{24}A\left( x\right)
\end{align*}%
where $A\left( x\right) \equiv 4x^{3}-9\pi x^{2}+\frac{9\pi ^{2}}{2}x$ and $%
B\left( x\right) \equiv x^{4}-3\pi x^{3}+\frac{9\pi ^{2}}{4}x^{2}.$ It is
easily verified that $A\left( 0\right) =A\left( \frac{3\pi }{2}\right) =0$
and $B\left( 0\right) =B\left( \frac{3\pi }{2}\right) =0$ so that%
\begin{align*}
H_{3}\left( 0\right) +E\left( 0\right) & =f\left( 0\right) \\
H_{3}\left( \frac{3\pi }{2}\right) +E\left( \frac{3\pi }{2}\right) &
=f\left( \frac{3\pi }{2}\right) \\
H_{3}^{\prime }\left( 0\right) +E^{\prime }\left( 0\right) & =f^{\prime
}\left( 0\right) \\
H_{3}^{\prime }\left( \frac{3\pi }{2}\right) +E^{\prime }\left( \frac{3\pi }{%
2}\right) & =f^{\prime }\left( \frac{3\pi }{2}\right)
\end{align*}%
Hence, $H_{3}\left( x\right) +E\left( x\right) $ has the \textit{same
interpolatory properties} (\ref{interp prop f}) and (\ref{interp prop f'})
as the original Hermite polynomial $H_{3}\left( x\right) .$

\subsection{Numerical integration}

Another obvious application is numerical integration, although we mention
this only briefly. With $E\left( x\right) $ approximated via a cubic spline,
we find%
\begin{align*}
\left\vert \int_{0}^{3\pi /2}f\left( x\right) dx-\int_{0}^{3\pi
/2}H_{3}\left( x\right) dx\right\vert & \sim 0.7 \\
\left\vert \int_{0}^{3\pi /2}f\left( x\right) dx-\int_{0}^{3\pi /2}\left(
H_{3}\left( x\right) +E\left( x\right) \right) dx\right\vert & \sim 3\times
10^{-12}
\end{align*}%
Clearly, there is a significant difference in accuracy and, of course, since 
$H_{3}\left( x\right) $ and $E\left( x\right) $ are polynomials, their
integrals are determined exactly.

\section{Conclusion}

We have shown how the Rolle function in Hermite interpolatory polynomial
approximation can be determined by solving an appropriate initial-value
problem. Consequently, the approximation error can be determined. In
particular, once the Rolle function is known, the Rolle term in the
expression for the approximation error can itself be approximated by means
of a polynomial, and this can result in a significant improvement in the
quality of the Hermite approximation overall. We have demonstrated this
effect using both a least-squares fit and a cubic spline, and we have
observed improvements in the accuracy of the approximation of many orders of
magnitude. This speaks to the potential value of the idea presented here,
and in \cite{Prentice 1}. We have also briefly observed that subsequent
numerical integration can also be made substantially more accurate, although
we will reserve further developments in that regard for future research.

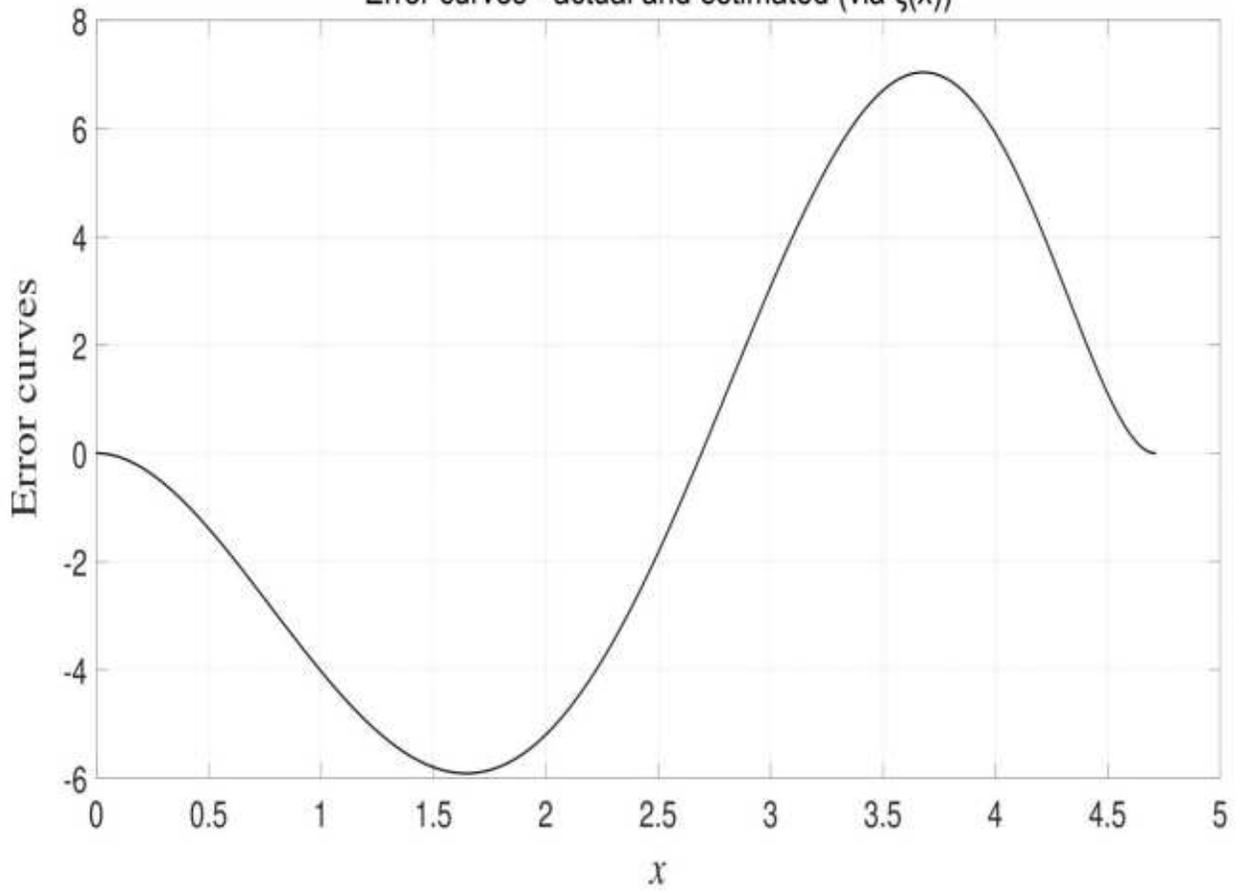
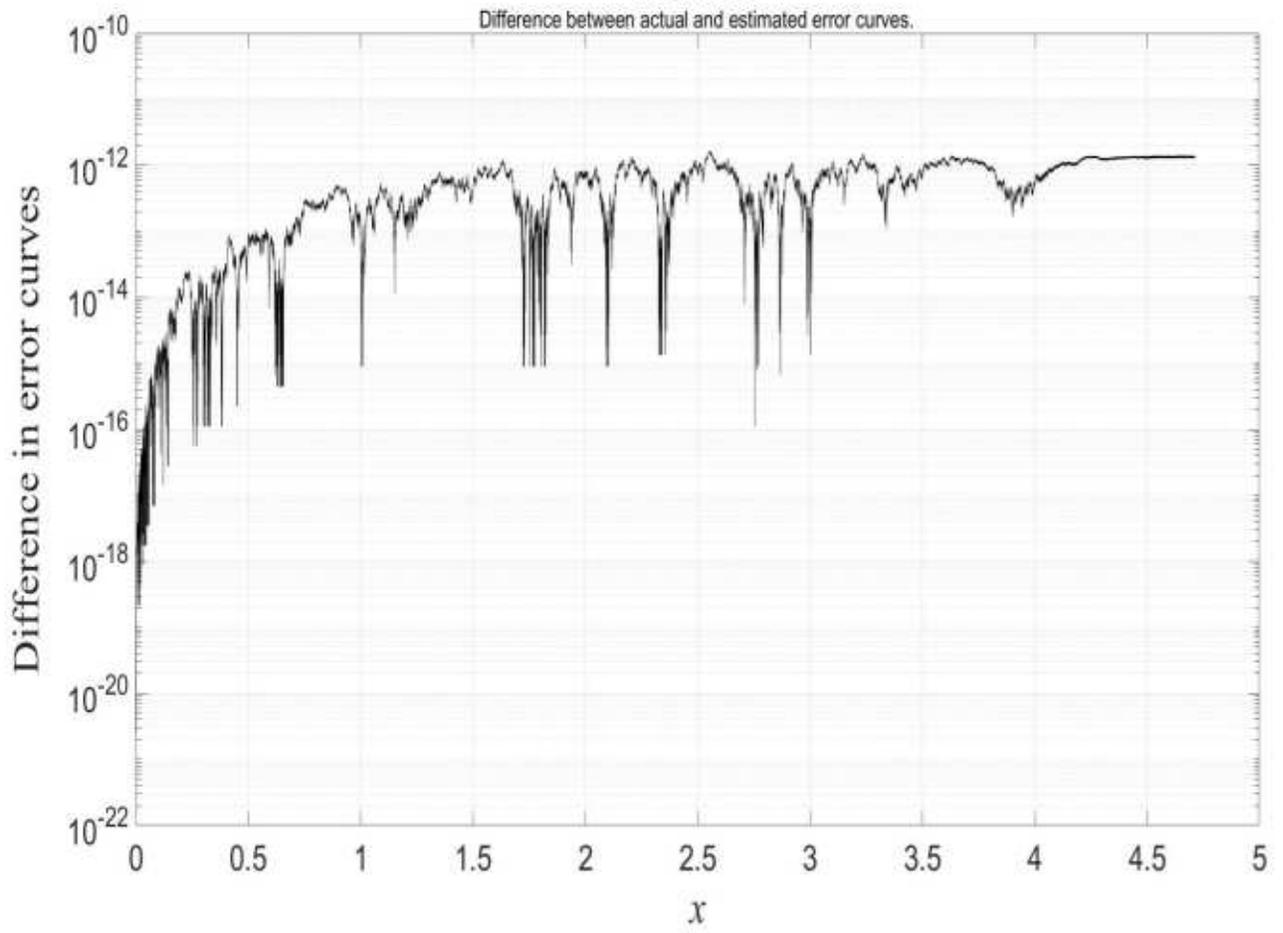
\end{document}